\pdfoutput=1

\input mypapers.mac
\input auxi.mac

\def\lf[#1]{{\lfloor #1 \rfloor}}
\def\blf[#1]{\bigg\lfloor #1 \bigg\rfloor }

\outputformat{dvi}{}  
\input drivers.mac

\documentclass{hebtex}

\begin{document}

\norevhebuse 
\nolabelnote

\twelverom   
\automathbaselineskip 

\forcecolortrue 



\begintitle 
{A Framework for Asymptotic Limit Problems of Probabilistic Nature}
\endtitle

\author{Michael Bensimhoun}{}{}{\date, Jerusalem}



\beginabstract

A convenient framework for dealing with asymptotic limit problems of 
probabilistic nature is provided.
These problems include questions such as finding the asymptotic proportion of terms 
of a sequence falling inside a given interval, or the limit of the
arithmetic mean of its partial sums; but several classes of problems 
are examined in a much more general setting.  
The proposed framework, which aims to unify those questions and their solution,
is based on the idea that to any finite multiset §E_n§,
one can associate a finitely distributed atomic probability §\mu_n§; 
assuming §\mu_n§ tends in distribution to a probability §\mu§,
it provides the tools needed to establish the desired asymptotic limit. 
Few examples are worked out in order to illustrate how using the framework.
\endabstract

\medbreak

 

\begincontent
\item{\sref{intro}.} Introduction 
\item{\sref{notations}.} Notations and Conventions
\item{\sref{multisets}.}  Multisets as Finitely Distributed Atomic Probabilities
\item{\sref{StieltjesNd}.} Stieltjes Integral in one or Several Dimensions 
\item{\sref{portmanteau}.}  Random Vectors and the Portmanteau Theorem 
\item{\sref{UsingFramework}.} Completing the Framework and Examples
\endcontent




\section intro introduction 

The idea of applying the theory of probabilities to compute certain limits is of course 
not new, nor even recent. 
Equidistribution theory for example, is more than one hundred years old, and 
has been examined and developed under the probabilistic framework 
(see for example the introduction of \cite{Li}).

Nevertheless, when it comes to solve certain classes of asymptotic limit problems of
probabilistic flavor,  
it may not be so easy for the non specialist to pick the suitable
results from the extensive literature, and to chain them together 
in order to obtain a formal proof. 
In fact, while most mathematicians would probably find the correct route for 
solving questions of this kind, a unifying and general exposition may still be 
of some value, if only to teach students some nice application of the 
theory of probability.    

Here are examples of questions of the aforementioned sort that 
are inquired further, where we 
denote by §\{x\}§ the fractional part of a real number §x§.
\beginlist[\bullet]

\item 
{\sl To find}
$$
\lim_{n\to \infty} {1\over n} \sum_{k = 1}^n \sin  \left(\Bigset{\sqrt k}\right);
$$

\item
{\sl 
to determine the asymptotic proportion of terms of the sequence
$$
u_k = \sin  \left(2\pi\Bigset{\sqrt k}\right)
$$
falling in §[-1/2, 1/2]§};

\item 
{\sl 
to find, for every §x\in [0,1]§, the asymptotic proportion of elements 
of the form §\{n/i\}§ that fall inside 
the interval §[0, x]§, as §n\to \infty§;
}

\item 
{\sl to find 
$$
\lim _{n\to \infty} {1\over n} \sum_{i=1}^n \gset {{n\over i}},
$$
and more generally} 
$$
\lim_{n\to \infty} {1\over n} \sum_{i=1}^n f\left( \gset{{n\over i}} \right);
$$

\item 
{\sl Investigate problems of the form
$$
\lim_{n\to \infty} {\displaystyle{1\over h(n)}} \sum_{i = 1}^{N(n)} f \big(P(i)/n\big),
$$
where §P§ and §h§ are polynomials, §f§ is measurable,
§N(n)§ is the greatest integer §i§ such that §P(i)\leq n§,
and the leading coefficients of §P§ and §h§ are §>0§.
}

\endlist

The framework presented here is composed of three parts; the first one deals with 
finitely distributed measures arising from multisets, the second part is  
an invocation of the multi-dimensional Stieltjes integral, and the 
third part is related to random vectors and the portmanteau theorem. 
The last section completes and explains how to use the framework for solving 
problems of the aforementioned kind. 

As a modest additional bonus (Example~3 below), 
we exhibit a formula in terms of the 
Digamma function, for the 
asymptotic proportion of terms of the fractional part of §n/i§ 
falling inside a sub-interval of §[0,1]§, as §n§ tends to~§\infty§.  
This explicit formula does not seem to be well known, at the very least.%
\footn{One of the most precise study of these topics is done in \cite{SV}; the authors 
have left the expression of their result in the form of an analytical series (1.9).}


\section notations Notations and Conventions

If §A§ is a Borel set of §\bbR^k§, and §\mu§ is a measure on §A§, then §\mu§
can be extended to the whole Borel §\sigma§-algebra of §\bbR^k§ by defining 
§\mu(B) = 0§ for every Borel set §B§ disjoint from §A§. So, for our purpose, 
there is no point to consider measures on subsets of~§\bbR^k§.

Accordingly, we shall note §\Omega = \bbR^k§, and denote by §I_\Omega§ the 
identity map §\Omega\to \Omega§.

If §x§ and §y§ are vectors of §\bbR^k§, we shall also write §x\leq y§ 
to mean §x_i \leq y_i§ for every §i§, where §x_i§ and §y_i§ are the §i§-th component 
of §x§ and §y§ \resp.

Assume that $E$ is a finite multiset of §\bbR^n§. 
We use the notation 
$$
E = \multiset{ a, a, b, c, d, d, d}
$$
to denote the elements §a,b,c,d§ of §E§ with their multiplicities.

We also denote §x\in \in E§ to mean ``§x§ belongs to §E§, including repetitions'', while 
§x\in E§ should be understood in the usual set theory meaning.
For example, if §E= \multiset{ 1,1,1, 2,2, 3, 4}§, then 
$$
\multiset{k\st k\in\in E\and k\leq 2} = \multiset{1,1,1, 2,2}
$$
and 
$$
\multiset{k \st k\in E \and k\leq 2} = \multiset{1,2} = \set{1,2},
$$
where we identify a multiset whose elements are all of multiplicity §1§ to the set containing the same
elements.

In general, every set theory symbol like §\incl§ should be understood in the set theoretic meaning 
(\thatis is defined with the relation §\in§), 
while the same symbol repeated twice, like §\incl\incl§ should be understood in the multiset meaning. 
So, for example, 
$$
E \incl [0,1] \iif x\in [0, 1] \quad \forall x\in E.
$$

The multi-cardinality of §E§, denoted by §\norm{E}§, is the number of elements of §E§, including repetitions, 
while §\absv{E}§ denotes the cardinality of the set of elements of §E§.


\section multisets Multisets as Finitely Distributed Atomic Probabilities
 
Let $(E_n)_n$ be a sequence of finite sets, or even of finite multi-sets, 
with §E_n\subset \Omega§.
Since a set can be viewed as a multiset for which
the multiplicity of each element is equal to~§1§, we shall  
suppose that §E_n§ is a multiset. 

To each $E_n$ corresponds a Borel atomic probability $\mu_n$ on $\Omega$ defined by
$$
\mu_n = {1\over \norm{E_n}}\sum_{e\in \in E_n} \delta_e,
$$
where $\delta_e(A)$ is equal to §1§ if §e\in A§, and to §0§ otherwise.  
Thus
$$
\mu_n(A) = {1\over \norm{E_n}} \Bignorm{ \multiset{ e\in \in E_n\st e\in A} }.
$$
\smallbreak

To fix the ideas, suppose for example that §(u_k)_k§ is a sequence of real numbers 
into §\Omega = [0,1]§.
Then we could define 
$$
E_n = \multiset{u_k\st 1\leq k\leq n},
$$ 
and for every Borel set §A§ of §\Omega§, §\mu_n(A)§ would be the proportion of terms
of §(u_n)_n§ which belong to §A§, among all the terms whose indices are §\leq n§.

As another typical example, we could define 
$$
E_n = \set{i/n\st 1\leq i \leq n},
$$ 
and §E_n§ could be seen as a §n§-sampling of the uniform probability on §[0,1]§.
Then §\mu_n(I)§ would be an approximation of the length of §I§, of order §O(1/n)§, for 
every subinterval §I§ of §\Omega§.
\smallbreak

More generally, in place of the uniform weights §1/\norm{E_n}§ above, 
we could consider \emph{a priori} weighted measures of the form 
$$
\mu_n = \sum_{e\in \in E_n} \omega_e \delta_e, \With \omega_e\geq 0,\quad \sum_{e\in \in E_n} \omega_e = 1,
$$
and this will be our final setting, where it is understood that if the weights §\omega_e§ are
not specified, they are simply equal to §1 / \norm{E_n}§.



\section StieltjesNd Stieltjes Integral in one or Several Dimensions

The Stieltjes integral on the real line is well known and we shall not repeat the basic theory here. 
We only wish to recall some more or less well known facts that are not always salient in classical expositions.  

As a preliminary remark,
if §\varphi\from [a,b]\to \bbR§ is of bounded variation, then it can be extended to the
whole of §\bbR§ by defining 
$$
\tilde \varphi(x) =\cases{ \varphi(b), & §x\geq b§\cr \varphi(a), & §x\leq a§\cr}.
$$
Then the variation of §\tilde\varphi§ \fnote{The variation of §\tilde\varphi§ is defined as the supremum 
of the sums §\sum_{[\alpha_i,\beta_i]\in S} \absv{\varphi(\alpha_i)- \varphi(\beta_i)}§ over all the finite sets §S§ 
of non-overlapping intervals
§[\alpha_i, \beta_i]§ of §\bbR§.}
is equal to the variation of §\varphi§ in §[a,b]§.
Moreover, it follows from the definition of the Stieltjes integral 
of a real function §f§ with respect to §\varphi§ in §[a,b]§, 
that §\int_a^b f \, d\varphi = \int_{-\infty}^\infty f \, d\tilde\varphi§, where
§f§ is extended in any manner outside §[a,b]§.
If in addition, §f§ is supposed to be continuous in §[a,b]§, then we can always 
extend it to the whole of §\bbR§ is such a way it be continuous in §\bbR§, 
by defining it, for example, to be equal to §\lim_{x\to a+} f(x)§ for §x\leq a§ and 
§\lim_{x\to b-} f(x)§ for §x \geq b§.
 
For these reasons, there is no point to include bounded intervals §[a,b]§ in the discussion below, and 
we shall always consider the functions are defined in the whole of §\bbR§, extending them in a suitable
way if necessary.

Here are three salient theorems regarding the Stieltjes integral, that 
are useful in practice.

\beginlist[\bullet]

\item 
{\sl If §\varphi§ is of bounded variation in §\bbR§, §f§ is regulated on the 
extended line §\bbRbar§
and continuous at all but finitely many points where §\varphi§ is continuous, 
then §\int f d\varphi§ exists.}

\item 
{\sl Assume §f§ and §\varphi§ have limits at §\pm \infty§.
If §\int f\, d\varphi § exists, then §\int \varphi\, df§ exists automatically, and 
$$
\int f \, d\varphi = [f(x)\varphi(x)]_{-\infty}^\infty  - \int \varphi \, df 
$$
(integration by parts).
}

\item  {\sl Whenever §\varphi§ is of bounded variation in §\bbR§,
there exists a unique signed Lebesgue measure §\mu§ on the Borel §\sigma§-algebra of 
§\bbR§, for which §\mu([a,b]) = \varphi(b) - \varphi(a)§ (§a, b\in \bbR, a\leq b)§.
Then the Stieltjes integral §\int f \, d\varphi§ coincides with §\int f\, d\mu§, the 
usual (signed) integral of §f§ with respect to §\mu§. 
Thus, whenever §\Var\, \varphi < \infty§, this so-called Lebesgue-Stieltjes integral
extends the Stieltjes integral and enriches it with the whole mechanics of 
measure theory. 
For a full generalization though, the Kurzweil-Henstock-Stieltjes 
integral is needed.} 

\item 
{\sl If $\varphi$ is continuously differentiable in §[a,b]§, then 
$$
\int_a^b f(t)\, d\varphi(t) = \int_a^b f(t)\, \varphi'(t)\, dt
$$
if one of these integrals exist.
Whenever §\varphi§ is absolutely continuous, this formula remains valid if the integrals
are interpreted as Lebesgue-Stieltjes integral\fnote{In fact, if $\varphi$ is only differentiable, 
the formula is valid in the context of the Kurzweil-Henstock-Stieltjes integral.}.
}
\endlist

\remark The integration by part formula above is not valid for the Lebesgue-Stieltjes integral.
Nevertheless, there exists versions of this formula for this integral as well, like the following statement (\cite{Re}).
\medbreak
{\sl Assume that both §f§ and §\varphi§ are of bounded variation (in particular they have limits at §\pm \infty§). 
Then if one of the following Lebesgue-Stieltjes integral
exists, the second integral exists as well, and there holds}
$$
\int f d\varphi = [f(x)\varphi(x)]_{-\infty}^\infty - \int \varphi_- \, df, \Where \varphi_-(t) = \lim_{s\to t-} \varphi(s).
$$


\subsection StieltjesNd Multidimensional Stieltjes integral

The Stieltjes integral has also been defined in the multi-dimensional context, a fact that does not seem to be widely known;
here, we only give some idea in the \hbox{2-dimensional} case, referring to \cite{Me} for a thorough exposition of the general case. 

In order to understand the idea behind the multi-dimensional Stieltjes integral,
suppose first that §F(x,y)§ is a two dimensional primitive of a 
continuous real function §f§ in §\bbR^2§:
$$
F(x,y) = \int_a^x \int_b^y f(x,y) \,dx\,dy.
$$  
It turns out that the integral of §f§ in the rectangle
$$
\calR = \set{x_1 < x \leq x_2,\quad y_1 < y \leq y_2}
$$
is equal to 
$$
\Delta_{\calR} F = F(x_2, y_2) - F(x_2, y_1) - F(x_1, y_2) + F(x_1, y_1),
$$
as can be seen easily with some picturing.

This fact generalizes to §n§ dimensions: the integral of §f§ in a hyper-rectangle 
obtains as a certain alternating sum of the n-dimensional primitive of §f§ at the vertices of 
the rectangle.

Now, noting §h_x = (x_2 - x_1)§ and §h_y = (y_2-y_1)§ and using Taylor's theorem, observe that
$$
\openup 2\jot
\eqalign{
\Delta_{\calR} F &= h_y {\partial F(x_2,y_1)\over \partial y} - h_y {\partial F(x_1, y_1)\over \partial y} 
		+ h_y^2 g(x_2, y_1) - h_y^2 g(x_1, y_1)\cr
&= h_xh_y {\partial^2 F(x_1,y_1)\over \partial x \partial y} + h_xh_y o(h_x, h_y)\cr
& = h_x h_y f(x_1, y_1) + h_x h_y \, o(h_x, h_y).  
}
$$
Hence 
$$
{\Delta_\calR f \over h_xh_y } \to f(x_1, y_1).
$$
We see that if a 2-dimensional Stieltjes §\int g \,df§ integral were defined somehow as 
the limit of the Riemann-Stieltjes sums
$$
\sum g \Delta_\calR F,
$$ 
then we would have
$$
\int g \,dF = \int g f \, dx = \int g {\partial ^2 F\over \partial x\partial y} dx.
$$
This suggests to define, for any function §\varphi§, 
$$\Delta_{\calR} \varphi = \varphi(x_2, y_2) - \varphi(x_2, y_1) - \varphi(x_1, y_2) + \varphi(x_1, y_1), $$
and the Stieltjes integral of a function §f§ with respect to §\varphi§ as the limit of the Riemann sums whose 
generic term is of the form 
$$
f(\xi) \Delta_{\calR} \varphi, \quad \xi\in \calR.
$$
Then if §\varphi§ is twice continuously differentiable, it fulfills, for every continuous function~§f§, 
$$
\int f d\varphi = \int f {\partial^2 \varphi\over \partial x \partial y} dx\,dy.
$$
The above definition of the 2-dimensional Stieltjes integral and the last relation generalizes to §n§ dimensions,
\emph{mutatis mutandis}:
\medbreak

{\sl
If a function §\varphi\from \Omega \to \bbR§ is §\nu§-times continuously\fnote{Assuming only §\varphi§ differentiable, 
the formula is still valid for the Kurzweil-Henstock-Stieltjes integral, 
and gives rise a Stieltjes-to-Lebesgue formula. See \cite{He} for the theory of the Kurzweil-Henstock integral in 
several dimensions.} 
differentiable, there holds, for every 
real function~§f§, 
$$
\int f(x_1, \ldots x_n)\; d\varphi(x_1, \ldots x_n) 
= \int f(x_1, \ldots x_n) \; {\partial^n \varphi(x_1,\ldots x_n) \over \partial x_1 \cdots \partial x_n}\; dx_1 \ldots dx_n 
$$
whenever one of these integrals exist.
}
\medbreak


\subsection nDVar Multivariate Variation and multi-dimensional Lebesgue Stieltjes Integral

The concept of \emph{variation} of a function in §\bbR^n§ also extends along the lines evoked in the previous section.
In the 2-dimensional case, we define 
$$
\absv{ \Delta_{\calR} }\, \varphi = \absv{ \varphi(x_2, y_2) - \varphi(x_2, y_1) - \varphi(x_1, y_2) + \varphi(x_1, y_1) }.
$$ 
Then the variation of §\varphi§ is defined as the supremum of all 
the sums whose generic term is of the form 
$$
\absv{ \Delta_{\calR_i}}\, \varphi, 
$$
over all the possible finite sets of non-overlapping rectangle §\calR_i§ of §\Omega§.

Again, this definition generalizes to the n-dimensional case, replacing rectangles by hyper-rectangles, and 
the above expression by a suitable alternating sum.

\medbreak

The Lebesgue-Stieltjes integral also generalizes to §n§-dimensions. 
Assuming §\varphi§ is of bounded variation,
there exists\fnote{To this end, the n-dimensional Stieltjes integral can be used to show that 
§\int_A d\varphi§ is a measure
on the algebra generated by the hyper-rectangles, with §\int_{\calR} d\varphi = \Delta_\calR \varphi§.} 
a unique measure §\mu§ on the Borel §\sigma§-algebra of §\bbR^n§ such that
§\Delta_\calR \varphi = \mu(\calR)§. 
So, as in the one-dimensional case, we can define the Lebesgue-Stieltjes integral  
of §f§ with respect to §d\varphi§ to be §\int f\,d\mu§. 
It coincides with the usual Stieltjes integral whenever this later exists, and is noted 
§ \int f \,d\varphi§ as well.



\section portmanteau Random Vectors and the Portmanteau Theorem

We recall that a \emph{random vector} §X§
is a pair §(h, \mu)§, where §h§ is a function 
§\Omega \to \bbR^\nu§, and §\mu§ is a probability on §\Omega§. 
By abuse, §h§ and §X§ are noted with 
the same letter, and we say that \emph{§X§ is a random vector on the probability space §(\Omega, \mu)§}.

Denoting §X = (X_1, \ldots, X_\nu)§, the \emph{cumulative distribution} of §X§ is the function §\bbR^\nu\to \bbR_+§
defined by
$$
\varphi(x_1, \ldots x_\nu) = \mu \Big( \set { X_i \leq x_i  \mtext{for all §i§} } \Big).
$$
It is of bounded variation, as can be seen with ease thanks to the §\sigma§-additivity of §\mu§.

By definition, the \emph{cumulative distribution of the measure §\mu§} is the cumulative distribution of the random vector §(I_\Omega, \mu)§.  
Explicitly, it is the function
$$
x\mapsto \mu(A_x) \With A_x = (-\infty, x_1]\times (-\infty, x_2]\times \cdots (-\infty, x_k],\quad x = (x_1, \ldots x_k).
$$

The \emph{expectation} of §X§ is defined by 
$$
E \set{X} = \int_\Omega X(x_1, \ldots x_k) \, d\mu. 
$$ 

The Portmanteau theorem states that the following conditions are equivalent, and defines the notion 
of \emph{convergence in distribution} of sequence of random vectors §(X_n)_n§ to a random vector §X§.
\goodbreak

\proposition [(\cite{Va}, Lemma~2.2 and \cite{Du}, \thm 3.10.1)] portmanteau
\beginlist 
\increasekeywidth by 1pt
\iitem   The sequence §(X_n)_n§ converges in distribution to a random vector §X§ of §\bbR^\nu§;
\iitem   §E\set{f(X_n)} \to E\set{f(X)}§ for every bounded continuous function §f§, 
\iitem  §E\set{f(X_n)} \to E\set{f(X)}§ for every bounded measurable function §f§ whose set 
of discontinuities is §\mu§-null; 
\iitem   §E\set{ f(X_n) } \to E\set{ f(X) }§ for every bounded Lipschitz  function §f§; 
\iitem   §E \set{\exp (i\; t\cdot X_n)} \to E\set{\exp(i\;t\cdot X)}§ for every §t\in \bbR^\nu§; 
(L\'evy's continuity theorem, \cite{Du}, \thm 3.10.5 or \cite{Va}, \thm 2.13);
\iitem   §\limsup_{n\to\infty} E\set{f(X_n)} \leq E\set{f(X)}§ for every upper semi-continuous 
function §f\from \bbR^\nu\to \bbR§ bounded from above; 
\iitem   §\liminf_{n\to \infty} E\set{f(X_n)} \geq E\set{f(X)}§ for every lower semi-continuous function
 §f\from \bbR^\nu\to \bbR§ bounded from below; 
\iitem §\limsup _{n\to \infty} \mu(X_n\in K)\leq \mu(X\in K)§ for every closed set §K§;
\iitem §\liminf_{n\to \infty}\mu(X_n\in G)\geq \mu(X\in G)§ for every open set §G§;
\iitem   §\varphi_n (x) \to \varphi(x)§ at every point §x§ where §\varphi§ is continuous, §\varphi§ being the 
cumulative distribution of §X§, and §\varphi_n§ that of §X_n§ (\cite{Va}, Lemma~2.2.1 or \cite{Du}, \thm 3.10.2); 
\iitem for every Borel set §A§ such that §\mu(X \in \partial A) = 0§, §\mu(X_n\in A)\to \mu(X\in A)§; 
\endlist 
\endproclaim

In what follows, we say, as is usual, that §(\mu_n)_n§ \emph{converges in distribution} to §\mu§, if the associated sequence of 
random vectors §(I_\Omega, \mu_n)§ converge in distribution to §(I_\Omega, \mu)§. According to \iref{x}, this 
means that §\mu_n(A_x) \to \mu(A_x)§, with
$$
A_x = (-\infty, x_1]\times (-\infty, x_2]\times \cdots (-\infty, x_k],\quad x = (x_1, \ldots x_k),
$$
at every continuity point §x§ of the function §x\mapsto \mu(A_x)§.
\medbreak

Now, we assume that the measures §\mu_n§ are those arising from multisets, as in \sec \sref{multisets}, and 
we also assume the definitions and notations there.

Suppose that §X_n§ is a random vector on §(\Omega, \mu_n)§, with values in §\bbR^\nu§ (§\nu\in \bbN^*)§:
$$
X_n = (X_{n,1}, X_{n,2},\ldots X_{n,\nu}).
$$
With the notations of \sec \sref{multisets}, let 
$$
E'_n = \multiset{ X_n(e)\st e\in\in E_n } \subset \bbR^\nu.
$$
Notice that §\norm{E'_n} = \norm{E_n}§ for all §n§.

Then §E'_n§ determines an atomic Borel measure §\mu'_n§ of § \bbR^\nu §, of the form 
$$
\sum_{e' \in \in E'_n} \omega'_{e'} \delta_{e'}, \autoeqno[mu2mu']
$$
where §\omega'_{e'} = 1/\norm{E'_n}§ in the usual case where the weighting is uniform. 
In the general case, §\omega'_{e'}§ is just constraint to fulfill the condition 
$$
\sum_{y\in \in E'_n,\  y = e'} \omega'_y = \sum_{ x \in \in E_n ,\  X_n(x) = e'} \omega_x,
$$
for which \eqref{mu2mu'} does not depend on the choice of the weights.

Moreover, there holds, for every Borel set §A\in \bbR^\nu§, 
$$
\mu'_n(A) = \mu_n(X_n \in A).
$$ 
Indeed, 
$$\openup 2\jot
\eqalign{
\mu'_n (A) &=  \sum _{e'\in \in E'_n} \omega'_{e'}  1_A(e')  
= \sum_{e'\in \in E'_n,\;  e'\in A} \omega'_{e'}  \cr
& = \sum_{e'\in A}\ \sum_{e\in \in E_n,\;  X_n(e)=e'} \omega_e 
 = \sum_{e\in\in E_n,\; X_n(e) \in A} \omega_e\cr
&= \mu_n(X_n \in A),
}
$$ 
which shows our contention.

As a consequence, §\mu'_n§ is the push forward measure of §\mu_n§, and by the well known 
push forward measure theorem, it follows that for every measurable function §f§ from §\bbR^\nu§ to some space, 
the expectation of the random vector §f(X_n)§ is
$$
E\set{f(X_n)} = \int_{\bbR^\nu} f(x_1, \ldots x_\nu)\, d \mu'_n.
$$
\medbreak 

Using random variables in §\bbR§, or more generally random vectors in §\bbR^k§, is more natural, 
allows more notational flexibility, and provides the right framework for the kind of questions we are
dealing with in this paper.   
Random vector are flexible in the sense that their composition is also a random vector. 
So, assuming for example that §\Omega'§ is a domain containing all the §E'_n§ and that 
$Y_n$ is a random vector §(\Omega',\mu'_n) \to R^m§, it is plain that §Y_n\circ X_n§ is a random vector 
§\Omega \to R^m§, and 
$$
\int Y_n d\mu'_n = \int Y_n\circ X_n d\mu_n.
$$ 
\medbreak

Denote by §\phi_n§ the cumulative distribution of §X_n§.
Since §\phi_n§ is also the cumulative distribution of §(I_{\bbR^\nu}, \mu'_n)§ in §\bbR^\nu§,
and since §\mu'_n§ is finitely distributed, it is clear 
that §\phi_{n}§ is a finitely supported step function, increasing with 
respect to each variable separately, with\fnote{The second relation is almost trivial
in dimension §1§, but is a bit more difficult in several dimensions. The proof involves 
the fact that if the values of §\varphi§ at any two adjacent vertices of §\calR§ 
along the same dimension are equal, then 
§\Delta_{\calR} \phi = 0§. Hence the variation of §\phi§ occurs only at its ``corner'' points.
}  
$$
\phi(\infty,\ldots, \infty) = 1 \And \Var \phi_n = 1. 
$$

Moreover, it follows from the definition of the multidimensional Lebesgue-Stieltjes integral that  
for every measurable function §f§ from §\bbR^\nu § to some space,
$$
 \sum_{e\in \in E_n} \omega_{e} f(X_n(e))
= \int_\Omega f(X_n)\, d \mu_n 
= \int_{\bbR^\nu} f\, d \mu'_n 
= \int_{\bbR^\nu} f\, d\phi_n
= \sum_{e'\in \in E'_n} \omega'_{e'} f(e'). 
\autoeqno[chainRel]
$$
Thus, according to a remark above,
$$
E\set{ f(X_n) } = \int_ {\bbR^\nu} f(x_1, \ldots x_\nu)\, d\phi_n.
$$
With §f = I_{\bbR^\nu}§, it follows that the expectation of §X_n§ is 
$$
E\set{ X_n } = \int_{\bbR^\nu} I_{\bbR^\nu}\, d\phi_n. 
$$
\medbreak

We now rewrite the portmanteau theorem above, specializing it to 
§X_n = (I_\Omega, \mu_n)§ and §X = (I_\Omega, \mu)§. 
As above, we denote by §\varphi_n§  the 
cumulative distribution of §\mu_n§, and by §\varphi§ the cumulative distribution of §\mu§.
\goodbreak

\proposition portmanteau2
The following statements are equivalent.
\beginlist 
\increasekeywidth by 1pt
\iitem   The sequence §\mu_n§ converges in distribution to  §\mu§;
\iitem   §\int f \, d\mu_n \to \int f d\mu§ for every bounded continuous function §f§;
\iitem   §\int f\, d\mu_n \to \int f d\mu§ for every bounded measurable function §f§ whose
set of discontinuities is §\mu§-null; 
\iitem   §\int f\, d\mu_n \to \int f\, d\mu § for every bounded Lipschitz  function §f§; 
\iitem   §\int_x \exp (i\; t\cdot x)\, d\mu_n \to \int_x \exp(i\;t\cdot x)\, d\mu§ for every §t\in \Omega§; 
\iitem   §\limsup_{n\to\infty} \int f\, d\mu_n \leq \int f\, d\mu§ for every upper semi-continuous 
function \newline §{f\from \Omega \to \bbR}§ bounded from above; 
\iitem   §\liminf_{n\to \infty} \int f\, \mu_n \geq \int f\, d\mu§ for every lower semi-continuous function
\newline §{f\from \Omega \to \bbR}§ bounded from below; 
\iitem   §\limsup_{n\to\infty} \mu_n(K)\leq \mu(K)§ for every closed set §K\in \Omega§;
\iitem   §\liminf_{n\to\infty} \mu_n(G)\geq \mu(G)§ for every open set §G\in \Omega§;
\iitem   §\varphi_n (x) \to \varphi(x)§ at every point §x§ where §\varphi§ is continuous;
\iitem   §\mu_n(A)\to \mu(A)§ for every continuity set §A§ of §\mu§,
\thatis every Borel set §A§ for which §\mu(\partial A) = 0§.  
\endlist 
\endproclaim

For what we are dealing with, the most important fact is the equivalence of \iref{i}--\iref{v} with \iref{x} 
and \iref{xi}.

{\sc Remark:}
The convergence in distribution of §\mu_n§ to §\mu§ does not imply that §\mu_n(A)\to \mu(A)§ for every Borel set. 
For example, if the multisets §E_n§ are of the form 
$$\set{i/n\st i=1,2,\ldots n},$$
 then surely §\mu_n(\bbQ) = 1§
for every §n§. But §\mu_n§ tends in distribution to the Borel measure §\mu§ on §[0,1]§, for which, of course, §\bbQ§ 
is a negligible set.
 
In fact, even if §A§ is a countable union of disjoint intervals of §[0,1]§, §\mu_n(A)§ may not tend to §\mu(A)§!
To see this, it suffices to observe that the union of the sets §E_n§ above can be covered by countably 
many disjoint open intervals, whose union §U§ is of total Borel measure §\epsi < 1§. 
But §\mu_n(U) = 1§ for all §n§, showing that §\mu_n(U)§ does not tends to §\mu(U)§. 
This does not contradict \iref{xi} because the boundary of §U§ is easily seen to be §[0,1]§, 
certainly not a null set for §\mu§.

We can nevertheless formulate the following remark, which is useful in practice.

\proposition contSet 
For a continuous measure §\mu§ and a Borel set §A§,
if both the set of non-interior points of §A§ and the set of accumulation points of §A§ are 
countable, or more generally, are countable unions of null sets, then §A§ is a continuity set of §\mu§.
\endproclaim

\proof This is clear since the boundary of §A§ is the union of the two sets 
mentioned in the proposition. \endproof

The following corollary to \thm \proclref{portmanteau2} will be used in some questions below.

\corollary portmantCor 
Assume that §g§ is a measurable function §\Omega \to \bbR^\nu§.
Assume as above that §\mu_n§ converges to §\mu§ in distribution, and let 
§\phi_n§ and §\phi§ be the cumulative distributions of the push forward measures §\mu'_n§ and §\mu'§ with 
respect to §g§:
$$
\phi_n(t) = \mu_n\big( \set{ g(x) \leq t} \big) \And \phi(t) = \mu\big( \set{ g(x) \leq t}\big).
$$
If the sets 
$$
A_t = \set{ x \st g(x)\leq t}
$$
are continuity sets of §\mu§, then 
the random vectors §X_n = (g, \mu_n)§ converge in distribution to §X = (g, \mu)§,
and the measures §\mu'_n§ converge in distributions to §\mu'§.

In particular, for every continuous function §f§  from §\bbR^\nu§ to some space,
$$\int_{\bbR^\nu} f\, d\phi_n  \to \int_{\bbR^\nu} f d\phi.$$ 
\endproclaim

\proof
Since §A_t§ are continuity sets for §\mu§, the portmanteau theorem \iref[portmanteau2]{xi} ensures that 
$$
\mu_n(A_t) \to \mu(A_t) 
$$
for every §t§.
Consequently §\phi_n(t) \to \phi(t)§ for every §t§.
The portmanteau theorem again \iref[portmanteau]{x} implies that §X_n§ converges to §X§ in 
distribution. 
Moreover, since §\phi_n§ and §\phi§ are also the cumulative distributions of §\mu'_n§ and §\mu'§ \resp,
it follows that §\mu'_n§ converges to §\mu'§ in distribution.
The last assertion is a repetition of one assertion of the portmanteau theorem.
\endproof



\section UsingFramework Completing the Framework and Examples

We are almost done. But we wish to present one more theorem, which generalizes, \emph{mutatis mutandis},
to the n-dimensional case.
As it is not the most well known theorem, we provide a proof.

\theorem StieltjesCOnv
Denote §\Omega = \bbR§.
Suppose that §\varphi_n§ are functions of bounded variation §\Omega\to \bbR§, which converges to a 
function §\varphi§ pointwise. 
Assume that §\Var\, \varphi_n§ is bounded for every §n§. 
Then §\varphi§ is of bounded variation, and the sequence of variations §\Var\, \varphi_n§ converges to §\Var\, \varphi§. 
\endproclaim

\proof 
Consider an arbitrary finite set of non overlapping finite intervals §[\alpha_i, \beta_i]§. 
Then according to the hypothesis,
$$
\sum_i \absv{\varphi_n(\beta_i) - \varphi_n(\alpha_i)} \leq M, 
$$
for some bound §M§.
Since the set §\set{\alpha_i, \beta_i}§ is finite, for §n§ sufficiently large, the above sum 
will be as close as we please to the sum 
$$
\sum_i  \absv{\varphi(\beta_i) - \varphi(\alpha)}.
$$
Therefore this later sum is not larger than §M§, showing that §\varphi§ is of bounded variation.

Now, set §V = \Var \varphi§, and let §\epsi > 0§.
It is possible to find a finite set §S_1§ of intervals of the above form, such that 
$$
\Bigabsv{\sum_i \absv{\varphi(\beta_i)-\varphi(\alpha_i)} - V} < \epsi. \autoeqno[sumApprox]
$$
Furthermore, the above argument shows that for §n§ sufficiently large, 
$$
\Bigabsv{\sum_i \absv{ \varphi_n(\beta_i) - \varphi_n(\alpha_i)} -  \sum_i \absv{\varphi(\beta_i)-\varphi(\alpha_i)} } < \epsi. 
$$
Combining these facts shows that for every §n§ sufficiently large, 
$$ 
\Bigabsv{\sum_i \absv{\varphi_n(\beta_i)-\varphi_n(\alpha_i)} - V} < 2\epsi. \autoeqno[sumApprox2]
$$
Consequently, denoting by §V'§ the supremum of §\set{\Var\, \varphi_n}§, we have 
$$
V' \geq V.
$$ 
On the other hand, if §V'§ were §> V§, there would exist, for every §N\in \bbN^*§, 
a number §n\geq N§, and a finite set §S_2§ of non-overlapping intervals §[\alpha_{n,i}, \beta_{n,i}]§ 
such that 
$$
\sum_i \absv{\varphi_n(\beta_{n,i}) - \varphi_n(\alpha_{n,i}) } > V + {V'-V\over 2}.
$$
By mixing the intervals of §S_1§ and §S_2§, we can produce a finite set §S§ of non-overlapping intervals §[\alpha'_j, \beta'_j]§ 
covering §S_1\union S_2§, and which is 
finer than §S_1\union S_2§ in the sense that every interval of §S§ is
 a sub-interval of an interval of §S_1\union S_2§.

Consequently, 
$$
\sum_j  \absv{\varphi_n(\beta'_{j}) - \varphi_n(\alpha'_{j}) } > V + {V'-V\over 2}.
$$
This contradicts \eqref{sumApprox2}, whenever  §\epsi§ is chosen to be 
§< (V'-V) /2§, and §N§ (hence also §n§) is sufficiently large.  
Therefore §V' = V§.

Finally,  assertion \eqref{sumApprox2} again shows that §\Var \varphi_n§ tends to §V§,
completing the proof. 
\endproof

The relevance of this theorem for our concern is such:
with the notations of \sec \sref{portmanteau}, if the sequence 
of cumulative distributions §(\varphi_n)_n§ for the probabilities §\mu_n§ converges
to a function §\varphi§, then §\Var_\Omega \varphi = 1§ since §\Var_\Omega \varphi_n = 1§ for all §n§.
Consequently, there exists, according to the general measure theory,
a unique probability measure §\mu§ on §\Omega§, the so-called Lebesgue-Stieltjes measure, 
for which §\mu([u, v]) = \varphi(v)-\varphi(u)§, or, in several dimensions, §\mu(\calR) = \Delta_{\calR} \varphi§. 
Moreover, §\mu_n§ converges obviously to §\mu§ in distribution.
Equivalently, we have now random vectors §X_n = (I_\Omega, \mu_n)§ which converge in distribution to §X = (I_\Omega, \mu)§,
and this allows using the mechanics convergence in distribution.
In other words, the sole pointwise convergence of the functions §\varphi_n§ to §\varphi§ is sufficient to place us in the frame
of the portmanteau theorem above.
 
\medbreak

The previous observations and theorems are all what is needed to solve 
many limit problems, as we shall now show.

\subsection kindPb Problem of the form §\lim_{n\to \infty}{1\over \norm{E_n}}\sum _{e\in \in E_n} f(e)§, §f§ continuous

Given finite multisets §E_n§, it is asked to find 
$$
\lim_{n\to \infty}{1\over \norm{E_n}}\sum _{e\in \in E_n} f(e).
$$

We have to determine, if possible, 
$$
\varphi_{n}(x) = \mu_n ( I_\Omega \leq x), \And \varphi(x) = \lim_{n\to\infty}  \varphi_n(x).
$$
Once this is done, \thm \proclref{StieltjesCOnv} ensures that §\mu_n§ converges in distribution 
to a probability §\mu§, whose cumulative distribution is §\varphi§.
By the portmanteau theorem (\proclref{portmanteau2}), it follows that the desired limit is 
$$
\int f \, d\mu = \int_\Omega f(t)\, d\varphi(t),
$$
which can usually be computed using Stieltjes or Lebesgue-Stieltjes integration theorems.

Sometimes, it is difficult or impossible to compute directly the above integral. 
Another method may be used. 
Seeing the function §f§ as a random vector §(f, \mu)§, one can try to 
compute the cumulative distribution function of the push forward measure 
§ \mu'(A) = \mu(f^{-1}(A))§. 
It is given by 
$$
\varphi(t) = \mu( f \leq t) =\int_\Omega f^{-1}\big( \set{ x\st f(x) \leq t} \big) d\varphi(x).  
$$  
Then by the push forward measure theorem applied to the function §t\mapsto t§ of §\bbR^\nu§, we have 
$$
\int_\Omega f\, d\mu = \int_{\bbR^\nu} t \, d\mu' = \int_{\bbR^\nu} t\, d\varphi(t).
$$
This is the well known push forward integration technique, which has actually nothing to do with the previous exposition, but is 
still useful for the kind of questions we are dealing with here. 

If furthermore §\nu=1§, a further simplification occurs by integrating by parts:
$$
\int_\Omega f\, d\mu = \int t\, d\varphi(t) = [t\varphi(t)]_{-\infty}^{\infty} - \int \varphi(t) \, dt.
$$

The paradigm we have just described includes the following typical categories:

\beginlist[\bullet]
\item 
{\sl Problems of the form
$$
\lim\limits_{n\to \infty} {1\over n} \sum_{k=1}^n f(u_k),
$$
where §(u_k)_k§ is a sequence of real numbers and §f§ is continuous.
}
\medbreak

The relevant multisets are 
$$
E_n = \multiset{ u_k\st k\leq n};
$$

\item 
{\sl Problems of the form} 
$$
\lim\limits_{n\to \infty} {1\over n} \sum_{i = 1}^n f( i/n) 
\Orthat \lim\limits_{n\to \infty} {1\over n} \sum_{i=1}^n f(n/i).
$$
\medbreak

Here the relevant multisets are 
$$
E_n = \multiset{{i\over n}\st 1\leq i\leq n}.
$$
With §I_\Omega§ and §\varphi_n§ defined as in the previous sections, 
it is clear that §\varphi_n(t) \to t§, the cumulative distribution function of 
the uniform probability in §[0, 1]§, which we denote §\mu§; 
hence the solution to the first problem is  
$$
E\set{f((I_\Omega, \mu_n))} \to E\set{f((I_\Omega, \mu))} = \int_0^1 f(t)\, dt. 
$$
The second problem reduces in fact to the first one by replacing §f(t)§ by §f(1/t)§, 
so its solution is
$$
\int_0^1  f(1/t)\, dt. 
$$ 

\item 
{\sl problems of the form
$$
\lim_{n\to \infty} {\displaystyle{1\over h(n)}} \sum_{i = 1}^{N(n)} f \big(P(i)/n\big),
$$
where §P§ and §h§ are polynomials and §f§ is continuous.
}
\medbreak

This is the generalization of the previous category of problems.
In the formula, §N(n)§ is the greatest integer §i§ such that §P(i)\leq n§.
It is also supposed that the leading coefficients of §P§ and §h§ are §>0§.
\newline 
The relevant multisets are 
$$
E_n = \multiset { P(i) / n \st 1\leq i\leq N(n)}.
$$
We have to compute 
$$
\varphi_n(x) = {1\over N(n)} \bigabsv{ \set { i \leq N(n)\st P(i)\leq nx} }.
$$
For §n§ large, the equation §P(u) = nx§ has a unique solution  (since 
the sign of §P'(u)§ is constant near §+\infty§).
Asymptotically, this implies that if the degree of §P§ is §q§ and its leading coefficient is §a\ (> 0)§, 
then the solution §u§ is of the order 
$$
\root q\of {nx\over a}.
$$
In particular, the order of §N(n)§ is asymptotically,
$$
\root q \of {n\over a}.
$$
It follows that for every §x \in [0,1]§,
$$
\varphi_n(x) \to \varphi(x) = \root q\of x, \As n\to \infty.
$$
For §x§ outside §[0,1]§, it is easy to see that 
$$
\varphi_n(x) \to \cases{ 0, & §x\leq 0§,\cr 1, & §x\geq 1§.}
$$ 
We deduce the following limit §L§:
$$
\lim_{n\to \infty} {\displaystyle{1\over N(n)}} \sum_{i = 1}^{N(n)} f \big(P(i)/n\big)
= \int f \,d\varphi = \int_0^1 f\, d\varphi = \int_0^1 f(x) \root q\of x \, dx.
$$
Moreover, if the degree of §h§ is §r§ and its leading coefficient is §b\ (>0)§,
§h(n)§ is asymptotically equal to §\root r \of {n\over b}§.
Then the solution to the question is
$$\openup 2\jot
{\root q\of {n\over a} \over \root r \of {n\over b}}\,L = 
\cases{+\infty, & $q > r$, \cr 0, &$q < r$, \cr 
\root q \of {b\over a} \displaystyle \int_0^1 f(x) \root q\of x \, dx, & $q = r$.} 
$$

\item 
{\sl 
Multidimensional versions of the previous problems.
}
\medbreak

Despite the examples here are 1-dimensional, the 
above framework is valid for problems in n-d as well. There are many interesting examples, but in
order to keep the length
of this article to a reasonable size, they will not been presented here. 

\endlist

\medbreak
\example[1]
{\sl We denote by §\{x\}§ the fractional part of a real number §x§.
To determine} 
$$
\lim_{n\to \infty} {1\over n} \sum_{k = 1}^n \sin  \left(\Bigset{\sqrt k}\right).
$$ 

We set §u_k = \bigset{\sqrt k} § and
$$
E_n = \multiset{u_k \st 1\leq k \leq n }\incl [0, 1].
$$

The main task is a classic exercise which consists in determining the pointwise limit §\varphi§ of 
$$\varphi_n(t) = {1\over n} \Bigabsv{ \set{ k\leq n \st u_k \leq t} }  .$$ 

If §t \leq 0§, §\varphi_n(t) = 0§ since §u_k\geq 0§, and if §t \geq 1§, §\varphi_n(t) = 1§ 
since §u_k\in [0, 1]§. Hence §\varphi(t) = 0§ for §t\leq 0§ and §\varphi(t) = 1§ for §t\geq 1§.
We assume from now on that §t\in [0,1]§.

For every §k§, there exists a unique §m\in \bbN§ and §0\leq x<1§ such that
$$
m\leq \sqrt k = m + x. 
$$
Equivalently,
$$
m^2 \leq k = (m+x)^2 = m^2 + 2mx + x^2.
$$
We have to count the number §S_n§ of integers §k\leq n§ for which 
the corresponding §x§ in the above expression fulfills §x \leq t§. 
  
Clearly, for a given §m§,
there are §2mt + t^2§ such numbers.
Since §m\leq \sqrt k \leq n§, §m§ must vary between 
§1§ and §\sqrt n§.
So the result is
$$\openup 2\jot
\eqalign{
S_n = \sum_{m =1}^{\floor{\sqrt n}} (2mt + t^2) &= 2t \sum_{m=1}^{\floor{\sqrt n}} m + t^2 \floor{\sqrt n}  \cr
&= t\Big(\floor{\sqrt n}^2 + \floor{\sqrt n}\Big) +  t^2 \floor{\sqrt n}  = tn + O(\sqrt n).
}
$$
Consequently, for every §0\leq t\leq 1§
$$
\varphi_n(t) = {1\over n} S_n \to t, \quad \mtext{as}\quad n \to \infty. 
$$
In other words,  §\mu_n§ converges in distribution
to the uniform probability law §\mu§ in §[0,1]§.

Then the portmanteau theorem ensures that the desired limit is 
$$
\int_\Omega \sin(t)\, d\varphi(t) = \int_0^1 \sin(t) \,dt = 1 - \cos(1). 
$$  


\subsection probKind2 Problems involving the cumulative distribution of §\mu'§

There are at least two circumstances where the cumulative distribution function of §\mu'§ is necessary.

In these questions, we suppose that the probability measures §\mu_n§ converge in distribution to 
a probability §\mu§ which is known or has been determined. 
We denote, as previously, by §\varphi_n§ and §\varphi§ 
their cumulative distribution \resp. 

\beginlist[\bullet]
\item {\sl Problems similar to the questions dealt in the previous section, except
§f\from \Omega\to \bbR^\nu§ is only supposed to be measurable but not continuous.}
\medbreak

In this case, we cannot use directly the portmanteau theorem to assert
that 
$$
\lim_{n\to \infty}\int f d\mu_n = \int f d\mu.
$$
Nevertheless, it may be possible to show that the sets 
$$
A_t = \set{ f(X) \leq t} 
$$
are continuity sets of §\mu§
(to this end, \prop \proclref{contSet} may be just fine).
Denote by §\mu'_n§ and §\mu'§ be the push forward measures of §\mu_n§ and §\mu§ with respect to §f§ \resp,
and by §\phi_n§ and §\phi§ the cumulative distributions of §\mu'_n§ and §\mu'§ \resp.
Then by \cor \proclref{portmantCor}, we can conclude that §\mu'_n§ converges to §\mu§ in distribution, 
§\phi_n\to \phi§, and 
$$
\int fd\mu_n \to \int f \, d\mu = \int f \, d\varphi.
$$
 
\item {\sl Problems where it is asked to find the asymptotic proportion
of terms §f(e)§ that fall inside an interval (or hyper-rectangle in several dimensions) §\calR§, 
with §e\in\in E_n§ and §f\from \Omega\to \bbR^\nu§ measurable.} 
 
\medbreak
We have the measure §\mu_n§ associated to the multiset §E_n§, and the measure §\mu'_n§
associated to the multiset
$$
E'_n = \multiset{ f(e) \st e\in\in E_n}.
$$
\newline
The problem is to find
$$
\lim_{n\to \infty} \mu'_n(\calR).  
$$
Denote by §\mu'§ the push forward measure of §\mu§ with respect to §f§, and 
by §\phi§ its cumulative distribution function.
Let §t_1, t_2,\ldots§ be the vertices of §\calR§.
It can often be checked that
$$
A_i = \set{x\st f(x) \leq t_i}
$$
are continuity sets of §\mu§, and this will be assumed here.
Then §\mu_n\to \mu§ pointwise on these sets. 
Equivalently, §\mu'_n(t \leq t_i)\to \mu'(t\leq t_i)§ for every §i§.  
According to \sec \sref{StieltjesNd}, this implies 
$$\mu'_n(\calR) \to \mu'(\calR) = \Delta_\calR \phi.$$
If for example §\nu = 2§ and the vertices §t_i§ are ordered from the bottom-left
clockwise, then 
$$
\Delta_\calR \phi = \mu(A_{3}) - \mu(A_{2}) - \mu(A_{4}) + \mu(A_{1}). 
$$ 

\item 
{\sl Same problem as above, but where the hyper-rectangle §\calR§ is replaced by a measurable set §B§ 
of §\bbR^\nu§.}
\medbreak

With the definitions of the previous problem, we have to find
$$
\lim_{n\to \infty} \mu'_n(B).
$$ 
As in the previous problem, we assume here that
$$
A_t = \set{x\st f(x) \leq t}
$$
are continuity sets of §\mu§, as it is often the case.
Then the portmanteau theorem implies
$$
\lim_{n\to\infty}\phi_n(t) = \lim_{n\to \infty} \mu'_n\big(\set{\tau \leq t}\big) \to \mu'(\tau\leq t) = \phi(t).  
$$
In other words, §\mu'_n§ tends to §\mu'§ in distribution.
\newline
Now, in order for §\mu'_n(B)§ to tend to §\mu'(B)§, we have to make sure that §B§ is a continuity set 
for §\mu'§.
This will hold if one of these two conditions are fulfilled:
$$
\int _{\partial B} d\phi = 0, \Orthat \mu(f^{-1}(\partial B)) = 0.
$$
Then 
$$
\mu'_n(B) \to \mu'(B) = \int_B d\phi .
$$

\endlist

\medbreak
\example[2] 
{\sl Denoting by §\{x\}§ the fractional part of a real number §x§,
to determine the asymptotic proportion of terms of the sequence
$$
u_k = \sin  \left(2\pi\Bigset{\sqrt k}\right)
$$
falling in §[-1/2, 1/2]§.
}
\medbreak

Let us set 
$$
E_n = \multiset{\gset{\sqrt k} \st 1\leq k \leq n}.
$$
Example~1 above shows that the associated probabilities §\mu_n§ tend 
in distribution to the uniform probability law §\mu§ in §[0, 2\pi]§, whose cumulative distribution,
restricted to §[0, 2\pi]§, is §\varphi = {t \over 2\pi}§.

Thus, according to the explanations above, we have to 
find the cumulative distribution §\phi§ of $\sin(X)$ for $X$ uniformly distributed in $[0,2\pi]$.

Suppose first that §0\leq t \leq 1§.
Then, restricting §x§ to §[0, 2\pi]§,  
$$
\sin x\leq t \iif x\in [0, \arcsin t] \Orthat x \in [\pi - \arcsin t, 2\pi].
$$

If on the other hand §-1\leq t < 0§ (hence §x > \pi§), then 
$$
\sin x \leq t \iif x\in [\pi - \arcsin t , 2\pi + \arcsin t].
$$
As a result, we have 
$$
\phi(t) = \mu\big( \set{\sin x \leq t,\ x\in [0, 2\pi]} \big) = {\arcsin t\over \pi} + {1\over 2}. 
$$

Moreover, in both case, §f^{-1}\big((-\infty, t]\big)§ is a continuity set of §\mu§, because §\mu§ is continuous
and the boundary of this set is finite. 
Therefore the desired asymptotic proportion is 
$$
\phi(1/2) - \phi(-1/2) = {\arcsin (1/2)\over \pi} - {\arcsin(-1/2)\over \pi} = {1\over 3}.
$$

\medbreak
\example[3]
{\sl With §\{x\}§ denoting the fractional part of a real number §x§, 
to find, for every §x\in [0,1]§, the asymptotic proportion of elements 
of the form §\{n/i\}§ that fall inside 
the interval §[0, x]§, as §n\to \infty§.
}
\medbreak

We consider 
$$
E_n = \multiset{{i\over n}\st 1\leq i\leq n}\subset [0, 1] \And f(x) = \gset{{1\over x}}, \ (0,1]\to [0,1).
$$
Then 
$$
\gset{{n\over i}} = f\left({i\over n}\right).
$$
Obviously,  §\mu_n§ tends in distribution to the uniform probability law in §[0,1]§, and to §0§ 
outside this interval. This define a continuous measure §\mu§ in §\Omega§.

With §\mu'_n§ the push forward measure of §\mu_n§ with respect to §f§, and §\phi_n§ its cumulative distribution,
we have to find 
$$
\lim_{n\to \infty} \mu'_n(\tau \leq t) = \phi_n(t).
$$

The function §f§ is continuous in §(0,1]§, but does not have a left limit at §0§.
So, it cannot be extended to a continuous function in §\Omega§, and the problem belongs 
to the above category. We extend §f§ by 
§f(x) = 0§ for every §x\nin (0, 1]§, which makes by the way §f§ continuous at §1§. 

According to the above discussion, we have first to determine §\phi(t) = \mu(f \leq t)§.
We have
$$
f(x) \leq t \iif x \leq 0 \Orthat x > 1 \Orthat  m\leq {1\over x} \leq m+t, \with m\in \bbN^* \and t\in [0,1).
$$
The last condition is equivalent to 
$$
{1\over m+t} \leq  x \leq {1\over m}.
$$
Hence, denoting §A = (-\infty, 0] \union (1, +\infty)§, we have for every §t\in [0, 1)§, and in 
fact for every §t\in [0,1]§,
$$
\openup 2\jot
f^{-1}((-\infty,t]) =A \union \Union_{m\in \bbN^*} \left[ {1\over m+t},
{1\over m}\right] \And \phi(t) = \sum_{m\in \bbN} \left( {1\over m} - {1\over m+t}\right)\!. 
\autoeqno[solmainProb]
$$
Of course, if §t < 0§,
$$
f^{-1}((-\infty, t]) = \emptyset \And \phi(t) = 0  
$$
and if §t \geq 1§,
$$
f^{-1}((-\infty, t]) = \Omega  \And \phi(t) = 1,
$$
so §\phi§ is continuous in §\bbR§ and so is §\mu'§.

We now observe that for every §t§, 
the set of non interior points of §f^{-1}((-\infty, t])§ is countable, while its set of accumulation points is empty. 
By \prop \proclref{contSet}, it follows that it is a continuity set of §\mu§ for every §t§.
Therefore we have
$$
\phi_n(t) \to \phi(t)
$$
for every §t§.

The problem is solved by series \eqref{solmainProb} that we now express in analytic form.

We have, for every §t\in [0,1]§, 
$$
\phi(t) = \sum_{m\in \bbN^*}{t\over m(m+t)}.
$$
This function is differentiable, with 
$$
\phi'(t) = \sum_{m\in \bbN^*} {1\over (m+t)^2}.
$$
It turns out that §\phi'§ can be written
$$\phi'(x) = \zeta(2, x) - {1\over x^2}$$ for all $x>0$, 
where $\zeta(s,z)$ is the so-called \beginhyperlink{https://en.wikipedia.org/wiki/Hurwitz_zeta_function}\emph{Hurwitz Zeta function}\endhyper.
This last function has an integral representation, which allows us to write 
$$\phi'(x) = \int_0^\infty {t e^{-tx} \over 1-e^{-t}} dt  
- {1\over x^2},\quad (x > 0). \autoeqno[HurwitzZetaInt]$$ 

On the other hand, if the real part of $z$ is positive then the digamma function has the following integral 
representation due to Gauss:
$$\psi(z) = \int_0^\infty \left(\frac{e^{-t}}{t} - \frac{e^{-zt}}{1-e^{-t}}\right)\,dt.$$ 
Therefore 
$$\psi'(z) = \int_0^\infty t\frac{e^{-zt}}{1-e^{-t}}\,dt.$$
Thus, from \eqref{HurwitzZetaInt}, 
$$\phi'(x) = \psi'(x) - {1\over x^2}.$$

We deduce that 
$$\phi(x) = \psi(x) + {1\over x} + C.$$
The constant $C$ has to be determined with the condition
$$\phi(1) = 1 = \psi(1) + {1\over 1} + C = -\gamma + 1 + C,$$
where $\gamma$ is the Euler-Mascheroni constant.
It follows that $C = \gamma$, hence 
$$\phi(x) = \psi(x) + {1\over x} + \gamma.$$
Also, from the known Laurent expansion of the digamma function $\psi$, it 
follows that 
$$
\phi(x)=-\sum_{k\geq 1}\zeta(k+1)(-x)^k .
$$


\medbreak
\example[4]
{\sl To find 
$$
\lim _{n\to \infty} {1\over n} \sum_{i=1}^n \gset {{n\over i}},
$$
§\{x\}§ denoting the fractional part of §x§, and more generally} 
$$
\lim_{n\to \infty} {1\over n} \sum_{i=1}^n f\left( \gset{{n\over i}} \right).
$$

We have to determine §\phi_n§ and §\phi§, and to prove that §\phi_n\to \phi§.
The work has already been done in Example~3 above. 
So, we can assert that the solution to the problem is 
$$
\int_\Omega f \,d\phi = \int_0^1 f\, d\phi = \int_0^1 f(t) \phi'(t) \, dt 
= \int_0^1 f(t) \left(\psi'(t) - {1\over t^2}\right) dt,
$$
with §\psi§ the digamma function.
The function §\psi§ has the following known series representation, for all §t > 0§:
$$
\psi(t) = -\gamma + \sum_{n=0}^\infty \left(\frac{1}{n + 1} - \frac{1}{n + t}\right)
 = -\gamma + {t-1\over t} + \sum_{n=1}^\infty {t-1\over (n+1)(n+t)}. 
$$
Whenever §f(t) = t§, we deduce that the previous integral is equal to
$$
\eqalign{
\int_0^1 t\left( \psi'(t) - {1\over t^2} \right) dt &= \big[t\psi(t) + 1\big]_0^1 - \int_0^1 \left ( \psi(t) + {1\over t}\right) dt\cr
&= -\gamma + 1 + \big[ \ln(\Gamma(t))+\ln t \big]_0^1 = -\gamma  + 1,
}
$$
where the last equality follows from the functional equation 
$$
\ln\big(\Gamma(z)\big) = \ln\big(\Gamma(z+1)\big) - \ln z, 
$$
valid for all §z > 0§, and the fact that §\Gamma(2) = \Gamma(1) = 1§.

As a bonus, with §H_n§ the §n§-th harmonic number, we can compute 
$$
\openup 2\jot
\eqalign{
{1\over n}\sum_{1\leq i\leq n} \floor{{n\over i}} &=  H_n - {1\over n}\sum_{1\leq i \leq n} \gset{{n\over i}}\cr 
& = \ln n + \gamma  + O(n^{-1})  - (-\gamma + 1 + o(1)) \cr
&= \ln n + 2\gamma - 1 + o(1). 
}
$$

This is a weak version of a theorem of Dirichlet, which asserts that the remainder 
§o(1)§ is in fact an §O(n^{-{1\over 2}})§.



\acknowledgements

The author wishes to thank X-Rui from the Mathematics StackExchange, for having 
indicated how the probabilistic method could be used to solve the question 
dealt in Example~3 above; this was the starting point of this article. 
Regarding this question, he also guessed the expression of the solution in terms 
of the Digamma function, and sketched a proof different from ours  
(see \beginhyperlink{https://math.stackexchange.com/questions/4831482/distribution-of-the-fractional-parts-of-n-1-n-2-ldots-n-n-as-n-tends-to-i} this thread). \endhyper



\def\gobbleit#1{}%

\bigbreak\bigskip\centerline{\bf Bibliography}
\medskip \nobreak
\frenchspacing
\tenrom
\baselineskip = 12\scpt

\bibliographystyle{plain}

\bibliography{lim}

\begin{thebibliography}{1}

\bibitem{Du}
R.~Durrett and R.~Durrett.
\newblock {\em Probability: Theory and Examples}.
\newblock Cambridge Series in Statistical and Probabilistic Mathematics.
  Cambridge University Press, 2019.

\bibitem{Li}
Vlada Limic and Ned\v zad Limi\'c.
\newblock Equidistribution, uniform distribution: a probabilist's perspective.
\newblock {\em Probability Surveys}, 15:~131--155, 2018.

\bibitem{Me}
Richard~M. Meyer.
\newblock {\em n-Dimensional Riemann-Stieltjes Integral}.
\newblock Springer New York, New York, NY, 1979.

\bibitem{He}
Henstock Ralph.
\newblock {\em {L}ectures on the {T}heory of {I}ntegration}.
\newblock World Scientific, Singapore, 1988.

\bibitem{Re}
D.~Revuz and M.~Yor.
\newblock {\em Continuous Martingales and Brownian Motion}, volume 293,
  chapter~0.
\newblock Springer-Verlag, Berlin, 2004.

\bibitem{SV}
Bahman Saffari and R.~C. Vaughan.
\newblock On the fractional parts of §x/n§ and related sequences. {II}.
\newblock {\em Annales the l'institut Fourier}, 27,~\numero 2:~1--30, 1977.

\bibitem{Va}
A.W. van~der Vaart.
\newblock {\em Asymptotic Statistics}.
\newblock Asymptotic Statistics. Cambridge University Press, 2000.

\end{thebibliography}


\end{document}